\documentclass[11pt]{article}
\usepackage{amsfonts}
\usepackage[dvips]{graphicx}
\usepackage{amssymb,color}
\usepackage[latin1]{inputenc}
\usepackage{epic}
\usepackage[latin1]{inputenc}
\usepackage{enumerate}
\usepackage{color}
\usepackage{lettrine}
\usepackage{calc,pifont}
\usepackage[noindentafter,calcwidth]{titlesec}
\usepackage{amsmath}
\usepackage{amssymb}
\usepackage{amsthm}
\usepackage{subfigure}
\usepackage{lineno} 

\def\BBox{\kern  -0.2cm\hbox{\vrule width 0.2cm height 0.2cm}}

\newtheorem{theorem}{Theorem}[section]

\newtheorem{corollary}{Corollary}[section]
\newtheorem{proposition}{Proposition}[section]
\newtheorem{remark}{Remark}[section]

\title{Covering and 2-degree-packing numbers in graphs}

\author{Carlos A. Alfaro \footnotemark[1] \and
Christian Rubio-Montiel \footnotemark[2] \and
Adri{\' a}n V\'azquez-\'Avila  \footnotemark[4]}
\date{}

\begin{document}
\maketitle

\def\thefootnote{\fnsymbol{footnote}}
\footnotetext[1]{Banco de México, Mexico City, Mexico, {\tt alfaromontufar@gmail.com, carlos.alfaro@banxico.org.mx}.}
\footnotetext[2]{Divisi{\' o}n de Matem{\' a}ticas e Ingenier{\' i}a, FES Acatl{\' a}n, UNAM, {\tt okrubio@ciencias.unam.mx}.} 
\footnotetext[4]{Subdirecci{\' o}n de Ingenier{\' i}a y Posgrado, UNAQ, Querétaro City, Mexico, {\tt adrian.vazquez@unaq.edu.mx}.}

\begin{abstract}
In this paper, we give a relationship between the covering number of a simple graph $G$, $\beta(G)$, and a new parameter associated to $G$ which is called 2-degree-packing number of $G$, $\nu_2(G)$. We prove that$$\lceil \nu_{2}(G)/2\rceil\leq\beta(G)\leq\nu_2(G)-1,$$ for any connected simple graph $G$, with $|E(G)|>\nu_2(G)$, and we give a characterization of simple connected graphs which attains the inequalities.
\end{abstract}





\section{Introduction}\label{sec:intro}

In this paper, we consider finite undirected simple graphs. For any undefined terms see \cite{Bondy}. Let $G$ be a graph, we call $V(G)$ the vertex set of $G$ and denote by $E(G)$ the edge set of $G$. For a subset $A\subseteq V(G)\cup E(G)$, $G[A]$ denotes the subgraph of $G$ which is induced by $A$. The distance between two vertices $u$ and $v$ in a graph $G$ is the number $d_G(u,v)$ of edges in any shortest $v-u$ path in $G$ that joins $u$ and $v$; if $u$ and $v$ are not joined in $G$, then $d_G(u,v)=\infty$. The \emph{neighborhood} of a vertex $u\in V(G)$, denoted by $N_G(u)$, is a subset of $V(G)$ adjacent to $u$ in $G$. The set of edges incident to $u \in V(G)$ is denoted by $\mathcal{L}_u$. Hence, the \emph{degree} of $u$, denoted by $deg(u)$, is $deg(x)=|\mathcal{L}_u|$. The minimum and maximum degree of a graph $G$ is denoted by $\delta(G)$ and $\Delta(G)$, respectively. Let $H$ be a subgraph of $G$. The \emph{restricted degree} of a vertex $u\in V(H)$, denoted by $deg_H(u)$, is defined as $deg_H(u)=|\mathcal{L}_u\cap E(H)|$. 

An \emph{independent set} of a graph $G$ is a subset $I\subseteq V(G)$ such that any two vertices of $I$ are not adjacent. The \emph{independence number} of $G$, denoted by $\alpha(G)$, is the maximum order of an independent set.A \emph{vertex cover} of a graph $G$ is a subset $T\subseteq V(G)$ such that all edges of $G$ has at least one end in $T$. The \emph{covering number} of $G$, denoted by $\beta(G)$, is the minimum order of a vertex cover of $G$. This invariant is well known and intensively studied in a more general context and with different names, see for example \cite{AK06,AK06_2,AKMM01,AGAL,Huicochea,Tancer,MS11}. On the other hand, a \emph{k-degree-packing set} of a graph $G$ ($k\leq\Delta(G)$), is a subset $R\subseteq E(G)$ such that $\Delta(G[R])\leq k$. The \emph{k-degree-packing number} of $G$, denoted by $\nu_k(G)$, is the maximum order of a $k$-degree-packing set. We are interested when $k=2$, since $k=1$ is the \emph{matching number} of a graph.

The 2-degree-packing number is studied in \cite{CGCA,AGAL,Avila,Avila2} on a more general context, but with a different name, as 2-packing number. It is important to say that the definition of 2-packing in graphs has different meaning: A set $X\subseteq V(G)$ is called a 2-packing if $d_G(u,v)>2$ for any different vertices $u$ and $v$ of $X$, that is, the2-packing is a subset $X\subseteq V(G)$ in which all the vertices are in distance at least 3 from each other, see for example \cite{TOPP1991229}. Therefore, we call 2-degree-packing instead of 2-packing just in case of graphs.

In \cite{AGAL}, was proved for any simple graph $G$ it satisfies

\begin{equation}\label{des:inf}
\lceil \nu_{2}(G)/2\rceil\leq\beta(G).
\end{equation}
In this paper, we prove that for any simple graph $G$, with $|E(G)|>\nu_2(G)$, it satisfies
\begin{equation}\label{des:sup}
\beta(G)\leq\nu_2(G)-1.
\end{equation}
Hence, by Equations (\ref{des:inf}) and (\ref{des:sup}), we have the following:
\begin{theorem}\label{thm:intro}
	If $G$ is a simple connected graph with $|E(G)|>\nu_2(G)$, then $$\lceil \nu_{2}(G)/2\rceil \leq\beta(G)\leq\nu_2(G)-1.$$	
\end{theorem}

The main result of this paper is give a characterization of simple connected graphs that attain the upper and lower inequality of the Theorem \ref{thm:intro}.

\section{Some results}\label{sec:results}
In the remainder of this note, for the terminology, notation and missing basic definitions related to graphs, the reader may consult \cite{Bondy}. Only connected graphs with $|E(G)|>\nu_2(G)$ are  considered, since $|E(G)|=\nu_2(G)$ if and only if $\Delta(G)\leq2$. Moreover, we assume $\nu_2(G)\geq4$, since in \cite{AGAL} was proved the following:

\begin{proposition}\emph{\textbf{\cite{AGAL}}}\label{prop:helly}
Let $G$ be a simple connected graph with $|E(G)|>\nu_2(G)$, then $\nu_2(G)=2$ if and only if $\beta(G)=1$.
\end{proposition}

\begin{proposition}\emph{\textbf{\cite{AGAL}}}\label{prop:3,2}
	Let $G$ be a simple connected graph with $|E(G)|>\nu_2(G)$. If $\nu_2(G)=3$, then $\beta(G)=2$.
\end{proposition}

\begin{figure}[t]
	\begin{center}
		\subfigure[]{\includegraphics[height =2.5cm]{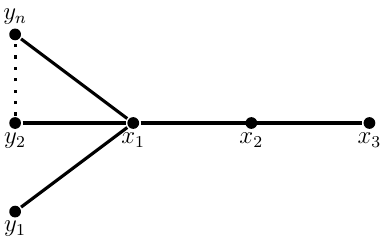}}\quad\quad\quad
		\subfigure[]{\includegraphics[height=2.5cm]{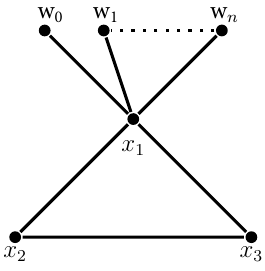}}
	\end{center}
	\caption{Graphs with $\beta=2$ and $\nu_2=3$.}
	\label{fig:nu_2_3}
\end{figure}

If $G$ satisfies the hypothesis of Proposition \ref{prop:helly}, then $G$ is the complete bipartite graph $K_{1,m}$. If $G$ satisfies the hypothesis of Proposition \ref{prop:3,2}, then $G$ is one of the graphs shown in Figure \ref{fig:nu_2_3} (see \cite{AGAL}).

The following Proposition \ref{remark:2-packing} show simple consequences of the definitions presented before, and some results are well known. 

\begin{proposition}\label{remark:2-packing}
	
	\	
	
	\begin{enumerate}
		\item If $R$ is a maximum 2-degree-packing of a graph $G$, then the components of $G[R]$ are either cycles or paths.
		
		\item If $G$ is either a cycle or a path, both of even length, and $T$ is a minimum vertex cover of $G$, then $T$ is an independent set.
		
		\item If $G$ is cycle of length odd and $T$ is a minimum vertex cover of $G$, then there exists an unique $u\in T$ such that $T\setminus\{u\}$ is an independent set. On the other hand, if $G$ is a path of length odd, then either there exists an unique $u\in T$ such that $T\setminus\{u\}$ is an independent set or $T$ is an independent and $\deg_T(u)=1$.	
		
		\item If $G$ is either a path or a cycle of length $k$, then $\beta(G)=\lceil \frac{k}{2}\rceil$.
		
		\item $\beta(K_n)=\nu_2(K_n)-1$.
	\end{enumerate}	
\end{proposition} 

\begin{remark}\label{remark}
Let $R$ be a maximum 2-degree-packing of a simple connected graph $G$. It is clear the number of components of $G[R]$ is at most $\nu_2(G)-1$. Moreover, if $T$ is a minimum vertex cover of $G[R]$, then $\beta(G)\leq k+p$, where $k$ is the number of components of $G[R]$ of a single edge, and $p=|\{v\in V(G[R]): deg_R(v)=2\}|$. Hence, $\beta(G)\leq k+p\leq\nu_2(G)$.
\end{remark}

\begin{proposition}
If $G$ is a simple connected graph with $|E(G)|>\nu_2(G)$, then $\beta(G)\leq\nu_2(G)-1$. 
\end{proposition}

\begin{proof}
Using the remark \ref{remark}, we have that $\beta(G)\leq k+p\leq\nu_2(G)$. If $k\geq1$, then it is not complicate to see that $\beta(G)\leq\nu_2(G)-1$. On the other hand, if $k=0$, then any component of $G[R]$ is a cycle, since if $G[R]$ has a path (of length at least 2) as a component, then $\beta(G)\leq\nu_2(G)-1$. Hence $p=\nu_2(G)$. We assume $V(G[R])=V(G)$, otherwise if $u\in V(G)\setminus V(G[R])$ and $e_u=uv\in E(G)\setminus R$, where $v\in V(G[R])$, then the following set $(R\setminus \{e_u\})\cup\{e_v\}$, where $e_v\in R$ is incident to $v$, is a maximum 2-degree-packing of $G$ with a path as a component, which implies that $\beta(G)\leq\nu_2(G)-1$. Therefore $\{v\in V(G[R]): deg_R(v)=2\}\setminus\{u\}$, for any $u\in V(G[R])$, is a vertex cover of $G$, implying that $\beta(G)\leq\nu_2(G)-1$. 
\end{proof}

Hence, we have the following:

\begin{theorem}\label{thm:des}
	If $G$ is a simple connected graph with $|E(G)|>\nu_2(G)$, then $$\lceil \nu_{2}(G)/2\rceil \leq\beta(G)\leq\nu_2(G)-1.$$	
\end{theorem}

\section{Graphs with $\beta=\nu_2-1$}\label{sec:tau=nu-1}
\begin{figure}
	\begin{center}
		\includegraphics[height =2cm]{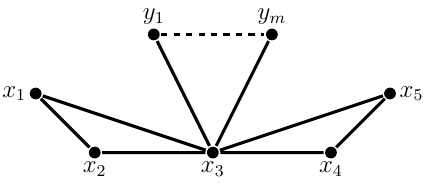}
		\hspace{1cm}
		\includegraphics[height=2.5cm]{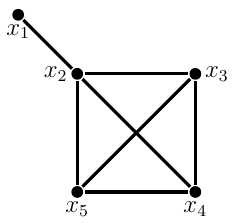}
	\end{center}
	\caption{Graphs with $\nu_2(G)=4$ and $\beta(G)=3$}\label{fig:nu_4}
\end{figure}

To begin with, some terminology is introduced in order to simplify the description of simple connected graphs $G$ such that $\beta(G)=\nu_2(G)-1$. 

In \cite{AGAL}, as a particular case, was proved the following:

\begin{proposition}\label{prop:nu=4}
If $G$ is a simple connected graph $G$ with $\nu_2(G)=4$ and $|E(G)|>4$, then $\beta(G)\leq3$.
\end{proposition} 

Moreover, in these same paper \cite{AGAL}, was given all simple connected graphs $G$ with $\nu_2(G)=4$ and $\beta(G)=3$, these graphs are certain subgraphs from Figure \ref{fig:nu_4} (see \cite{AGAL}). Hence, by Proposition \ref{prop:nu=4} we assume $\nu_2(G)\geq5$.

\begin{figure}[t]
	\begin{center}
		\subfigure[]{\includegraphics[height =2.5cm]{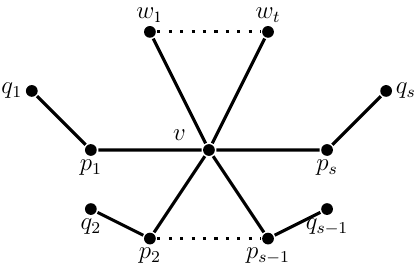}}
		\hspace{1cm}
		\subfigure[]{\includegraphics[height=2.5cm]{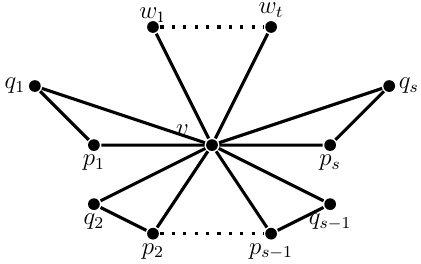}}
	\end{center}
	\caption{In $(a)$ depict the Graph $T^{r}_{s,t}$ and in $(b)$ depict the graph $G^r_{s,t}$.}
	\label{fig:arbol_gamma}
\end{figure}

In \cite{Avila3}, was defined the graph $T_{s,t}$, with $s\geq1$ and $t\geq2$, as follow (see Figure \ref{fig:arbol_gamma} $(a)$):
\begin{eqnarray*}
	V(T_{s,t})&=&\{p_1,\ldots,p_{s}\}\cup\{q_1,\ldots,q_{s}\}\cup\{w_1,\ldots,w_t\},\\
	E(T_{s,t})&=&\{p_iq_i:i=1,\ldots,s\}\cup\{vp_i:i=1,\ldots,s\}\cup\{vw_i:i=1,\ldots,t\}.
\end{eqnarray*}

And we define $G_{s,t}$, with $s\ge1$ and $t\geq2$, as follow  (see Figure \ref{fig:arbol_gamma} $(b)$):
\begin{eqnarray*}
	V(G_{s,t})&=&V(T_{s,t}),\\
	E(G_{s,t})&=&E(T_{s,t})\cup\{vq_i:i=1,\ldots,s\}.
\end{eqnarray*}

As a consequence of Corollary 2.4 of \cite{Avila3}:

\begin{corollary}\cite{Avila3}
$\beta(T_{s,t})=\nu_2(T_{s,t})-1=s+1$, for every $s\geq1$ and $t\geq2$.
\end{corollary}

Since the graph $T_{s,t}$ is a spanning graph of $G_{s,t}$, and any minimal vertex cover of $T_{s,t}$ is a vertex covering of $G_{s,t}$, we have the following:

\begin{corollary}
$\beta(G^{r}_{s,t})=\nu_2(G^{r}_{s,t})-1=s+1$, for every $s\geq1$ and $t\geq2$
\end{corollary}

\begin{corollary}
If $T_{s,t}$ is a spanning subgraph of a graph $G$ and $G$ is a spanning subgraph of $G_{s,t}$, then $\beta(G)=\nu_2(G)-1=s+1$.
\end{corollary} 

Let $R_1,\ldots,R_s,R_{s+1},\ldots,R_k$ be the components of a simple connected graph $G$, where $|R_i|=1$, for $i=1,\ldots,s$ and $|R_j|>1$, for $j=s+1,\ldots,k$. It is not difficult to see that $s\leq\nu_2(G)-2$. If $s=\nu_2(G)-2$, implies that $k=\nu_2(G)-1$ and $|E(G[R_k])|=2$. Hence, any edge from $E(G)\setminus E(G[R])$ is incident with the only one vertex $v\in V(G[R_k])$ with $deg_R(v)=2$. Hence, if $R_i=p_iq_i$, for $i=1,\ldots,s$, $R_k=w_0vw_1$, and $V(G)\setminus V(G[R])=\{w_3,\ldots,w_t\}$ (an independent set), if $t\geq3$, then $T_{s,t}$ is a spanning subgraph of a graph $G$ and $G$ is a spanning subgraph of $G_{s,t}$. Therefore, $\beta(G)=\nu_2(G)-1=s+1$. 

Let $R_1,\ldots,R_s,R_{s+1},\ldots,R_k$ be the components of a simple connected graph $G$, with $k$ as small as possible , where $|R_i|=1$, for $i=1,\ldots,s$ and $|R_j|>1$, for $j=s+1,\ldots,k$. Then, it is clear that $\beta(G)=s+\beta(H)$ and $\nu_2(G)=s+\nu_2(H)$, where $H$ is the graph defined as follow 
\begin{eqnarray*}
V(H)&=&\bigcup_{i=1}^kV(G[R_i])\cup \left(V(G)\setminus V(G[R])\right),\\
E(H)&=&E(G)\setminus\{R_1,\ldots,R_s\}.
\end{eqnarray*}
Hence, if $\tau(G)=\nu_2(G)-1$, then $\tau(H)=\nu_2(H)-1$. Therefore, we assume that any simple connected graph $G$, with $|E(G)|>\nu_2(G)$, has a maximum 2-degree-packing $R$ of $G$, where each component of $G[R]$ has at least 2 edges; and as consequence, the set $T=\{u\in V(G[R]):deg_{G[R]}(u)=2\}$ is a vertex cover of $G$.

Let $K_n^1$ be the simple connected graph defined as follow: 
\begin{eqnarray*}
V(K_n^1)&=&\{x_1,\ldots,x_n\}\cup\{u\},\\
E(K_n^1)&=&\{x_ix_j: 1\leq i<j\leq n\}\cup\{ux_1\}.
\end{eqnarray*}

The graph $K_n^1$ is the complete graph of $n$ vertices joined with an edge. It is easy to see that $\beta(K_n^1)=\nu_2(K_n^1)-1=n-1$.	

\begin{proposition}\label{thm:completa_graph}
Let $G$ be a simple connected graph with $|E(G)|>\nu_2(G)$, $\nu_2(G)\geq5$ and $\beta(G)=\nu_2(G)-1$. If $R$ is a maximum 2-degree-packing of $G$ with $V(G[R])=V(G)$, then either $G$ is the complete graph $K_{\nu_2}$ or $G$ is $K_{\nu_2}^1$, where $\nu_2=\nu_2(G)$.
\end{proposition}

\begin{proof}
Let $R$ be a maximum 2-degree-packing of $G$ with $V(G[R])=V(G)$, and let $R_1,\ldots,R_k$ be the components of $G[R]$, with $k$ as small as possible. Then
	\begin{itemize}
		\item [Case(i)] If $k=1$, then $G[R]$ is either a spanning path or a spanning cycle of the graph $G$. Let suppose that $R=u_0u_1\cdots u_{\nu_2-1}u_0$ is a spanning cycle: If there are two non-adjacent vertices $u_i,u_j\in V(G[R])$, then $T=V(G[R])\setminus\{u_i,u_j\}$ is a vertex cover of $G$ of cardinality $\nu_2(G)-2$, which is a contradiction. Therefore, any different pair of vertices of $G$ are adjacent. Hence, the graph $G$ is the complete graph of $\nu_2(G)$ vertices.
		
		On the other hand, if $R=u_0u_1\cdots u_{\nu_2}$ is a path, then $T=\{u_1,\ldots,u_{\nu_2-1}\}$ is a minimum vertex cover of $G$. Let assume that either $u_0u_j\in E(G)$ or $u_{\nu_2}u_j\in E(G)$, for all $u_j\in T^*=T\setminus\{u_1,u_{\nu_2-1}\}$, otherwise, $T\setminus\{u_j\}$ is a vertex cover of $G$ of cardinality $\nu_2(G)-2$, which is a contradiction. Without loss of generality, let suppose that $u_0u_j\in E(G)$, for all $u_j\in T^*=T\setminus\{u_1,u_{\nu_2-1}\}$. If $u_ju_{\nu_2}\in E(G)$, for some $u_j\in T^*$, then $R^*=(R\setminus\{u_ju_{j+1}\})\cup\left\{u_ju_{\nu_2}\right\}$ (since $\nu_2(G)\geq5$) is a 2-degree-packing with $G[R^*]$ as a cycle, which is a contradiction. Hence $u_ju_{\nu_2}\not\in E(G)$,	for all $u_j\in T^*$, which implies that $deg(u_{\nu_2})=1$. On the other hand, if there are two vertices $u_i\neq u_j\in T^*$ non-adjacent, then $\left(T\setminus\{u_i,u_j\}\right)\cup\{u_0\}$ is a vertex cover of $G$ of size $\nu_2(G)-2$, which is a contradiction. Also, $u_1u_j\in E(G)$ and $u_ju_{\nu_2-1}\in E(G)$, for all $u_j\in T^*$, otherwise there exists $u_j\in T^*$ such that either $(T\setminus\{u_1,u_j\})\cup\{u_0\}$ or $(T\setminus\{u_j,u_{\nu_2-1}\})\cup\{u_0\}$ is a vertex cover of $G$ of size $\nu_2(G)-2$, which is a contradiction. Therefore, the graphs $G$ is the graph $K_{\nu_2}^1$.
		
		\item [Case (ii)] Let suppose that $k\geq2$ and $T=\{v\in V(G[R]): degR(v)=2\}$. If there is at least one component as a paths (of length at least 2), say $R_1$, then
		\begin{eqnarray*}
			\beta(G)\leq|T|&\leq&(|E(R_1)|-2)+\sum_{i=2}^k|E(R_i)|\\
			&=&\sum_{i=1}^{k}|E(R_i)|-2=\nu_2(G)-2,
		\end{eqnarray*}
		which is a contradiction. Hence, $G[R_i]$ is a cycle, for all $i=1,\ldots,k$. 
		
		If there are two vertices $u,v\in V(G[R])$ such that $uv\not\in E(G)$, then $T\setminus\{u,v\}$ is a vertex cover of $G$ with $\beta(G)\leq\nu_2(G)-2$, which is a contradiction. Then, any two vertices $u,v\in V(G[R])$ are adjacent, which implies that $k=1$, a contradiction. Therefore, $G[R]$ is the complete graph of $\nu_2(G)$ vertices.   
	\end{itemize} 
\end{proof}

\begin{theorem}\label{coro:completa_pelo}
Let $G$ be a simple connected graph with $\nu_2(G)\geq5$ and $\beta(G)=\nu_2(G)-1$. Then either $G$ is the complete graph $K_{\nu_2}$ or $G$ is $K_{\nu_2}^1$, where $\nu_2=\nu_2(G)$.
\end{theorem}

\begin{proof}
Let $R$ be a maximum 2-degree-packing of $G$ and $I=V(G)\setminus V(G[R])$. Let assume that $I\neq\emptyset$, otherwise, the theorem holds by Proposition \ref{thm:completa_graph}. Hence, if $I\neq\emptyset$, then $I$ is an independent set of vertices. 
	\begin{itemize}
		\item [Case (i):] Let suppose that $G[R]$ is the complete graph of $\nu_2(G)$ vertices. We claimed that, if $u\in I$, then $deg(u)=1$. To verify the claim, let suppose on the contrary, $u$ is incident to at least two vertices of $V(G[R])$, say $v$ and $w$. If $V(G[R])=\{u_1,\ldots,u_{\nu_2}\}$, then without loss of gene\-rality, we suppose $u_1=v$ and $u_j=w$, for some $j\in\{2,\ldots,\nu_2\}$. Since $G[R]$ is a complete graph, then		
		$$(R\setminus\{u_1u_{\nu_2},u_{j-1}u_j\})\cup\{uu_1,uu_j,u_{j-1}u_{\nu_2}\}$$ is a 2-degree-packing of $G$ of size $\nu_2(G)+1$, which is a contradiction. Hence, if $u\in I$, then $deg_G(u)=1$. 
		
		On the other hand, if $|I|>1$, let $u,v\in I$. Without loss of generality, let suppose that $u$ is incident to $u_1$ and $v$ is incident to $u_j$, for some $j\in\{2,\ldots,\nu_2\}$. Since $G[R]$ is a complete graph, then	
		$$(R\setminus\{u_1u_{\nu_2},u_{j-1}u_j\})\cup\{uu_1,u_{j-1}u_{\nu_2},vu_j\}$$ is a 2-degree-packing of size $\nu_2(G)+1$, which is a contradiction. If $u$ and $v$ are adjacent to $u_1$, then $$(R\setminus\{u_1u_2,u_1u_{\nu_2}\})\cup\{uu_1,vu_1,u_2u_{\nu_2}\}$$is a 2-degree-packing of size $\nu_2(G)+1$, which is contradiction. Hence, $I=\{u\}$ with $deg(u)=1$, which implies that the graph $G$ is $K_{\nu_2}^1$.

		\item[Case (ii):] Let suppose that $G[R]$ is the graph $K_{\nu_2}^1$
		Let $v\in V(G)$ such that the $G[R]-v$ is the complete graph of size $\nu_2(G)$. If $u\in I$ is such that $uw\in E(G)$, whit $w\in V(G[R])$, then there exists a 2-degree-packing of $G$ of size $\nu_2(G)+1$ (see proof of Proposition \ref{thm:completa_graph}), which is a contradiction. Then $uw\not\in E(G)$, for all $w\in V(G[R])\cup\{v\}$, which implies that $G$ is a disconnected graph, unless $I=\emptyset$, and the theorem holds by Proposition \ref{thm:completa_graph}. 
	\end{itemize}
\end{proof}

\section{Graphs with $\beta=\displaystyle\left\lceil\nu_2/2\right\rceil$}\label{sec:tau_inf}
To begin with, some terminology and results are introduced in order to simplify the description of the simple connected graphs G which satisfy $\beta(G)=\lceil{\nu_2(G)/2\rceil}$.

\begin{proposition}\label{lema:ciclos_par}
Let $G$ be a simple connected graph and $R$ be a maximum 2-degree-packing of $G$. 
	\begin{enumerate} 
		\item If $\nu_2(G)$ is an even integer and $\beta(G)=\displaystyle\frac{\nu_2(G)}{2}$, then the components of $R$ has even length.
		\item If $\nu_2(G)$ is an odd integer and $\beta(G)=\displaystyle\frac{\nu_2(G)+1}{2}$, then there is an unique component of $R$ of odd length.
	\end{enumerate}  
\end{proposition}

\begin{proof}
We will prove 1., since the proof of 2. is completely analogous to  the proof of 1.: Let $R$ be a maximum 2-degree-packing of $G$ and let $R_1,\cdots,R_k$ be the components of $G[R]$. If $T$ is a minimum vertex cover of $G$, then $$\frac{\nu_2(G)}{2}=\beta(G)=|T|=\sum_{i=1}^{k}|T\cap V(R_i)|\geq\sum_{i=1}^{k}\beta(R_i)=\sum_{i=1}^{k}\lceil{\nu_2(R_i)/2\rceil}.$$
Hence, if $R_1$ have a odd number of edges, then$$\sum_{i=1}^{k}\lceil{\nu_2(R_i)/2\rceil}=\frac{\nu_2(R_1)+1}{2}+\sum_{i=2}^{k}\lceil{\nu_2(R_i)/2\rceil}\geq\frac{1}{2}+\sum_{i=1}^{k}\frac{\nu_2(R_i)}{2}=\frac{1}{2}+\frac{\nu_2(G)}{2},$$which is a contradiction. Therefore, each component of $G[R]$ has an even number of edges.
\end{proof}

Let $A$ and $B$ be two sets of vertices. The complete graph whose set of vertices is $A$ is denoted by $K_A$. The
graph whose set of vertices is $A\cup B$ and whose set of edges is $\{ab:a\in A, b\in B\}$ is denoted by $K_{A,B}$. On the other hand, let $k\geq3$ be a positive integer. The cycle of length $k$ and the path of length $k$ are denoted by $C^k$ and $P^k$, respectively.

If $A$ and $B$ are two sets of vertices from $V(C^k)$ and $V(P^k)$ (not necessarily disjoint) and $I$ be an independent set of vertices different from $V(C^k)$ and $V(P^k)$ then $C_{A,B,I}^k=(V(C_{A,B,I}^k),E(C_{A,B,I}^k))$ and  $P_{A,B,I}^k=(V(P_{A,B,I}^k),E(P_{A,B,I}^k))$ are denoted to be the graphs with $V(C_{A,B,I}^k)=V(C^k)\cup I$ and $V(P_{A,B,I}^k)=V(P^k)\cup I$, respectively, and $E(C_{A,B,I}^k)=E(C^k)\cup E(K_A)\cup E(K_{A,B})\cup E(K_{A,I})$ and $E(P_{A,B,I}^k)=E(P^k)\cup E(K_A)\cup E(K_{A,B})\cup E(K_{A,I})$, respectively. In an analogous way, we denote by $C_I^k$ to be the graph with $V(C_I^k)=V(C^k)\cup I$ and $E(C_I^k)=E(C^k)$ and we denote by $P_I^k$ to be the graph with $V(P_I^k)=V(P^k)\cup I$ and $E(P_I^k)=E(P^k)$. We define
$\mathcal{C}_{A,B,I}^k$ be the family of connected graphs $G$ such that $C_I^k$ is a subgraph of $G$ and $G$ is a subgraph of $C_{A,B,I}^k$. Similarly, we define $\mathcal{P}_{A,B,I}^k$ be the family of connected graphs $G$ such that $P_I^k$ is a subgraph of $G$ and $G$ is a subgraph of $P_{A,B,I}^k$. That is $$\mathcal{C}_{A,B,I}^k=\{G: C_I^k\subseteq G\subseteq C_{A,B,I}^k\mbox{ where $G$ is a connected graph}\}$$ $$\mathcal{P}_{A,B,I}^k=\{G: P_I^k\subseteq G\subseteq P_{A,B,I}^k\mbox{ where $G$ is a connected graph}\}$$

\begin{proposition}\label{prop:transversal_min_independence_set}
Let $k\geq4$ be an even integer, $T$ be a minimum vertex cover of $C^k$ and $I$ be an independent set of vertices different from $V(C^k)$. If $\hat{T}=V(C^k)\setminus T$ and	$G\in \mathcal{C}_{T,\hat{T},I}^k$, then $\beta(G)=\frac{k}{2}$ and $\nu_2(G)=k$.
\end{proposition}
\begin{proof}
It is clear that, if $G\in \mathcal{C}_{T,\hat{T},I}^k$, then $\beta(G)=\frac{k}{2}$. On the other hand, since $C^k$ is a 2-degree-packing of $G$, then $\nu_2(G)\geq k$. Moreover, since $\displaystyle\left\lceil\nu_2(G)/2\right\rceil\leq\beta(G)=\frac{k}{2}$, then $\nu_2(G)=k$.
\end{proof}

\begin{corollary}\label{coro:transversal_min_independence_set}
Let $k\geq4$ be an even integer, $T$ be a minimum vertex cover of $P^k$ and $I$ be an independent set of vertices different from $V(P^k)$. If $\hat{T}=V(P^k)\setminus T$ and $G\in \mathcal{P}_{T,\hat{T},I}^k$, then $\beta(G)=\frac{k}{2}$ and $\nu_2(G)=k$.
\end{corollary}

Now, let $\mathcal{\hat{C}}_{A,B,I}^k$ be the family of simple connected graphs $G$ with $\nu_2(G)=k$, such that $C_I^k$ is a subgraph of $G$ and $G$ is a subgraph of $C_{A,B,I}^k$. Similarly, let  $\mathcal{\hat{P}}_{A,B,I}^k$ be the family of simple connected graphs $G$ with $\nu_2(G)=k$ such that $P_I^k$ is a subgraph of $G$ and $G$ is a subgraph of $P_{A,B,I}^k$. That is  $$\mathcal{\hat{C}}_{A,B,I}^k=\{G: C_I^k\subseteq G\subseteq C_{A,B,I}^k\mbox{ where $G$ is connected and $\nu_2(G)=k$}\},$$ $$\mathcal{\hat{P}}_{A,B,I}^k=\{G: P_I^k\subseteq G\subseteq P_{A,B,I}^k\mbox{ where $G$ is connected and $\nu_2(G)=k$}\}.$$	

Hence if $k\geq4$ is an even integer, $T$ is a minimum vertex cover of either $C^k$ or $P^k$, and $I$ is an independent set different from either $V(C^k)$ or $V(P^k)$, then by Proposition \ref{prop:transversal_min_independence_set} and Corollary \ref{coro:transversal_min_independence_set} we have $$ \mathcal{\hat{C}}_{T,\hat{T},I}^k=\mathcal{C}_{T,\hat{T},I}^k\mbox{ and }\mathcal{\hat{P}}_{T,\hat{T},I}^k=\mathcal{P}_{T,\hat{T},I}^k.$$

However, if $k\geq5$ is an odd integer, $T$ is a minimum vertex cover of either $C^k$ or $P^k$ and $I$ is an independent set different from either $V(C^k)$ or $V(P^k)$, then $$ \mathcal{\hat{C}}_{T,\hat{T},I}^k\neq\mathcal{C}_{T,\hat{T},I}^k\mbox{ and }\mathcal{\hat{P}}_{T,\hat{T},I}^k\neq\mathcal{P}_{T,\hat{T},I}^k.$$ To see this, let $R$ be the cycle of length $k$ and $u,v\in T$ adjacent. Hence, if $G$ is such that $V(G)=V(C^k)\cup\{w\}$, where $w\in I$ and $E(G)=E(C^k)\cup\{uw,vw\}$, then $G\in \mathcal{C}^k_{T,\hat{T},I}$. However, it is clear that $\nu_2(G)=k+1$, which implies that $G\not\in\mathcal{\hat{C}}^k_{T,\hat{T},T}$. A similar argument is used to prove that $\mathcal{\hat{P}}_{T,\hat{T},I}^k\neq\mathcal{P}_{T,\hat{T},I}^k$.

\begin{proposition}\label{prop:ultimo}
Let $k\geq5$ be an odd integer, $T$ be a minimum vertex cover of $C^k$ and $I$ be an independent set of vertices different from $V(C^k)$. If $\hat{T}=V(C^k)\setminus T$ and	$G\in \mathcal{\hat{C}}_{T,\hat{T},I}^k$, then $\beta(G)=\frac{k+1}{2}$.
\end{proposition}

\begin{proof}
It is clear that

$$\frac{k+1}{2}=\displaystyle\left\lceil\nu_2(C^k_I)/2\right\rceil\leq\displaystyle\left\lceil\nu_2(G)/2\right\rceil\leq\beta(G)\leq|T|=\frac{k+1}{2},$$which implies that $\beta(G)=\frac{k+1}{2}$.
\end{proof}

\begin{corollary}\label{coro:ultimo}
Let $k\geq5$ be an odd integer, $T$ be a minimum vertex cover of $P^k$ and $I$ be an independent set of vertices different from $V(P^k)$. If $\hat{T}=V(P^k)\setminus T$ and $G\in \mathcal{\hat{P}}_{T,\hat{T},I}^k$, then $\beta(G)=\frac{k+1}{2}$.
\end{corollary}

\begin{proposition}\label{prop:caracterization}
Let $G$ be a connected graph with $|E(G)|>\nu_2(G)$ and $R_1,\ldots, R_k$ be the components of a maximum 2-degree-packing of $G$. If $\beta(G)=\left\lceil\nu_2(G)/2\right\rceil$, then $\beta(G)=\displaystyle\sum_{i=1}^k\beta(R_i)$.
\end{proposition}

\begin{proof}
Let $R$ be a maximum 2-degree-packing of $G$ and $R_1,\ldots,R_k$ be the components of $G[R]$. Since $R_i$ is a cycle or a path of length $\nu_2(R_i)$, then $\beta(R_i)=\left\lceil\nu_2(R_i)/2\right\rceil$, for $i=1,\ldots,k$. If $\beta(G)=\left\lceil\nu_2(G)/2\right\rceil$, then by Proposition \ref{lema:ciclos_par} we have $$\left\lceil\nu_2(G)/2\right\rceil=\beta(G)\geq\sum_{i=1}^k\beta(R_i)=\sum_{i=1}^k\left\lceil\nu_2(R_i)/2\right\rceil=\left\lceil\nu_2(G)/2\right\rceil.$$Therefore $\beta(G)=\displaystyle\sum_{i=1}^k\beta(R_i)$.
\end{proof}

By Proposition \ref{lema:ciclos_par} and Proposition \ref{prop:caracterization}, we have:

\begin{theorem}\label{coro:caracterization}
Let $G$ be a connected graph with $|E(G)|>\nu_2(G)$ and $R_1,\ldots, R_k$ be the components of a maximum 2-degree-packing of $G$. Then $\beta(G)=\left\lceil\nu_2(G)/2\right\rceil$, if and only if, $\beta(G)=\displaystyle\sum_{i=1}^k\beta(R_i)$, being 
\begin{enumerate}
	\item $|R_i|$ an even integer, for $i=1,\ldots,k$, if $\nu_2(G)$ an even number.
	\item $|R_1|$ is an odd integer and $|R_i|$ is an even integer, for $i=2,\ldots,k$, if $\nu_2(G)$ is an odd number.
\end{enumerate}
\end{theorem}

\begin{proposition}\label{prop:1}
Let $G$ be a simple connected graph with $\nu_2(G)\geq4$, $|E(G)|>\nu_2(G)$ and  $R_1,\ldots, R_k$ be the components of a maximum 2-degree-packing $R$ of $G$, with $k$ as small as possible. If $\beta(G)=\left\lceil\nu_2(G)/2\right\rceil$, then $I=I_1\cup\cdots\cup I_k=V(G)\setminus V(G[R])$, where either $I_i=\emptyset$ or for every $u\in I_i$ satisfies $N(u)\subseteq V(R_i)$, for $i=1,\ldots,k$.		
\end{proposition}
\begin{proof}
Let suppose that there exists $u\in I$, $w_i\in V(R_i)$ and $w_j\in V(R_j)$, for some $i\neq j\in\{1,\ldots,k\}$, such that $uw_i,uw_j\in E(G)$. Hence $(R\setminus\{e_{w_i},e_{w_j}\})\cup\{uw_i,uw_j\}$, where $w_i\in e_{w_i}\in E(R_i)$ and $w_j\in e_{w_j}\in E(R_j)$, is a maximum 2-degree-packing with less components than $R$, which is a contradiction. Therefore $I=I_1\cup\cdots\cup I_k$, where either $I_i=\emptyset$ or for every $u\in I_i$ satisfies $N(u)\subseteq V(R_i)$, for $i=1,\ldots,k$.		
\end{proof}

\begin{corollary}\label{coro:1}
Let $G$ be a simple connected graph with $\nu_2(G)\geq4$, $|E(G)|>\nu_2(G)$, $R_1,\ldots, R_k$ be the components of a maximum 2-degree-packing $R$ of $G$, with $k$ as small as possible, and $I=I_1\cup\cdots\cup I_k=V(G)\setminus V(G[R])$, where either $I_i=\emptyset$ or for every $u\in I_i$ satisfies $N(u)\subseteq V(R_i)$, for $i=1,\ldots,k$. If $\beta(G)=\left\lceil\nu_2(G)/2\right\rceil$, then $\beta(G[R_i])=\left\lceil\nu_2(G[R_i])/2\right\rceil$, for $i=1,\ldots,k$.			
\end{corollary}

\begin{proposition}\label{prop:2}
Let $G$ be a simple connected graph with $\nu_2(G)\geq4$, $|E(G)|>\nu_2(G)$ and $R$ be a maximum 2-degree-packingof $G$, such that $G[R]$ is a connected graph. If $\beta(G)=\left\lceil\nu_2(G)/2\right\rceil$, then either $G\in \mathcal{\hat{C}}_{T,\hat{T},I}^k$ or $G\in \mathcal{\hat{P}}_{T,\hat{T},I}^k$, where $T$ is a minimum vertex cover of either $C^k$ or $P^k$, $\hat{T}=V(G[R])\setminus T$ and $I=V(G)\setminus V(G[R])$.
\end{proposition}
\begin{proof}
By Proposition \ref{lema:ciclos_par}, we have either $\hat{C}^k_I$ is a subgraph of $G$ or $P^k_I$ is a subgraph of $G$. Let $T$ be a minimum vertex cover of $G$ (hence, a minimum vertex cover of $G[R]$, by Proposition \ref{prop:caracterization}). Hence, by definition, if $e\in E(G)\setminus E(G[R]$, then $e$ has an end in $T$, which implies that $G$ is a subgraph of $\hat{C}^k_{T,\hat{T},I}$. Therefore, either $G\in \mathcal{\hat{C}}_{T,\hat{T},I}^k$ or $G\in \mathcal{\hat{P}}_{T,\hat{T},I}^k$.  
\end{proof}

By Proposition \ref{prop:caracterization}, Proposition \ref{prop:2} and Corollary \ref{coro:1}, we have:

\begin{corollary}\label{coro:final}
Let $G$ be a simple connected graph with $\nu_2(G)\geq4$, $|E(G)|>\nu_2(G)$, $R_1,\ldots, R_k$ be the components of a maximum 2-degree-packing $R$ of $G$, with $k$ as small as possible, and $I=I_1\cup\cdots\cup I_k=V(G)\setminus V(G[R])$, where either $I_i=\emptyset$ or for every $u\in I_i$ satisfies $N(u)\subseteq V(R_i)$, for $i=1,\ldots,k$. If $\beta(G)=\left\lceil\nu_2(G)/2\right\rceil$, then either $G[V_i]\in \mathcal{\hat{C}}_{T_i,\hat{T}_i,I_i}^{k_i}$ or $G[V_i]\in \mathcal{\hat{P}}_{T_i,\hat{T}_i,I_i}^{k_i}$, where $V_i=V(G[R_i])\cup I_i$, $k_i=\nu_2(G[R_i])$, $T_i$ is a minimum vertex cover of either $C^{k_i}$ or $P^{k_i}$ and $\hat{T}_i=V(G[R_i])\setminus T_i$.			
\end{corollary}

Hence, by Proposition \ref{prop:transversal_min_independence_set}, Proposition \ref{prop:2}, Corollary \ref{coro:transversal_min_independence_set} and Corollary \ref{coro:final}, we have

\begin{theorem}
Let $G$ be a simple connected graph with $\nu_2(G)\geq4$, $|E(G)|>\nu_2(G)$, $R_1,\ldots, R_k$ be the components of a maximum 2-degree-packing $R$ of $G$, with $k$ as small as possible, and $I=I_1\cup\cdots\cup I_k=V(G)\setminus V(G[R])$, where either $I_i=\emptyset$ or for every $u\in I_i$ satisfies $N(u)\subseteq V(R_i)$, for $i=1,\ldots,k$. Then $\beta(G)=\left\lceil\nu_2(G)/2\right\rceil$, if and only if, either $G[V_i]\in \mathcal{\hat{C}}_{T_i,\hat{T}_i,I_i}^{k_i}$ or $G[V_i]\in \mathcal{\hat{P}}_{T_i,\hat{T}_i,I_i}^{k_i}$, where $V_i=V(G[R_i])\cup I_i$, $k_i=\nu_2(G[R_i])$, $T_i$ is a minimum vertex cover of either $C^{k_i}$ or $P^{k_i}$ and $\hat{T}_i=V(G[R_i])\setminus T_i$, being 
\begin{enumerate}
	\item $|R_i|$ an even integer, for $i=1,\ldots,k$, if $\nu_2(G)$ an even number.
	\item $|R_1|$ is an odd integer and $|R_i|$ is an even integer, for $i=2,\ldots,k$, if $\nu_2(G)$ is an odd number.
\end{enumerate}
\end{theorem}

{\bf Acknowledgment}

Research was partially supported by SNI and CONACyT.

\end{document}